%% file: manuscript.tex
\documentclass[final,onefignum,onetabnum]{siamart250211}

\input{ex_shared}

\ifpdf
\hypersetup{
  pdftitle={On the Superconvergence of ESFR Schemes},
  pdfauthor={Mathias Dufresne-Piché, and Siva Nadarajah}
}
\fi

\begin{document}

\maketitle

\begin{abstract}
  The energy stable flux reconstruction (ESFR) method provides an efficient and flexible framework to devise high-order linearly stable numerical schemes which can achieve high levels of accuracy on unstructured grids. While superconvergent properties of ESFR schemes have been observed in numerical experiments, no formal proof of this behavior has been reported in the literature. In this work, we attempt to address this by providing a simple derivation for the superconvergence of the dispersion-dissipation error of ESFR schemes for the linear advection problem when using an upwind numerical flux. We show that the superconvergence of ESFR schemes essentially relies on the capacity of the latter to generate superconvergent rational approximants of the exponential function, which is reminiscent of well-known theoretical results for superconvergence of discontinuous Galerkin (DG) methods. We also demonstrate that the drops in order of accuracy which are observed in numerical experiments as the ESFR scalar $c$ is increased are caused by both a modification of the structure of these rational approximants and a change in the multiplicity of the physical eigenvalue of the schemes as $c \to \infty$. Finally, our theoretical results are successfully validated against numerical experiments.
\end{abstract}



\section{Introduction}
Flux reconstruction (FR) schemes \cite{huynh2007flux, WANG20098161} are a class of high-order spatial discretizations which have been gaining in popularity in recent years as a simple and versatile finite-element formulation for computational fluid dynamics (CFD) problems. Various well-known schemes, such as collocation-based nodal discontinuous Galerkin (DG) methods and some spectral difference (SD) methods can be expressed within the FR unifying framework \cite{huynh2007flux}. Using the FR approach, Vincent et al.~\cite{vincent2011insights} were able to identify an infinite family of linearly stable methods known as energy stable flux reconstruction (ESFR) schemes, which have been the subject of extensive investigations in the past decade. The numerical properties of ESFR schemes for linear problems have been studied by \cite{vincent2011insights, lambert2023l2}, and ESFR discretizations have been developed for triangular and tetrahedral elements by \cite{williams2013energy, castonguay2012new} and extended to the advection-diffusion problem by \cite{quaegebeur2019stability1, quaegebeur2019stability2, quaegebeur2020insights}.

It is well-known that DG methods exhibit so-called superconvergent properties, which allow different error functionals to converge faster than the standard L2-error of the solution. For instance, \cite{adjerid2009discontinuous, yang2012analysis} showed that the DG spatial discretization of the linear advection problem was superconvergent with order $p+2$ at the Radau points. Superconvergence of order $2p+1$ was also observed by \cite{yang2012analysis} for the cell averages and the average error at the downwind points. A formal proof of this property was later provided by \cite{cao2014superconvergence}. Superconvergence of the dispersion-dissipation error of the DG method for linear hyperbolic problems has also been studied extensively. In \cite{hu2002eigensolution}, Hu and Atkins conjectured that the dispersion-dissipation relations of the DG spatial discretization were related to a specific Padé approximation of the exponential function, thereby leading to superconvergent dispersion error of order $2p+3$ and dissipation error of order $2p+2$. This conjecture was subsequently proven by \cite{ainsworth2004dispersive}, who investigated the different convergence regimes of the dispersion-dissipation error under $h$ and $p$ refinement. In \cite{guo2013superconvergence}, Guo et al. used von Neumann analysis to decompose the point-wise error associated with DG schemes. By symbolically computing the relevant eigenvalues and eigenvectors for low-order spatial discretizations, they were able to relate the superconvergent behavior of DG schemes to the spectral properties of the latter. This perspective was further investigated in the work of \cite{chalmers2017spatial}, which formally established the connection between a specific Padé approximation of the exponential function and superconvergence of the long-time error of DG schemes. \cite{chalmers2017spatial} also showed that this Padé approximation was related to the superconvergent behavior of the cell average and the error at the Radau and downwind points.

For FR schemes, superconvergence of the dispersion-dissipation error has been observed in numerical experiments by \cite{huynh2007flux, vincent2011insights, vermeire2016properties}. Namely, in \cite{vincent2011insights}, Vincent et al. measured that the order of accuracy of the dispersion-dissipation error of ESFR schemes was either $2p+1$, $2p$ or $2p-1$ depending on the value of the ESFR scalar parameter. The connection between the dispersion-dissipation error of FR schemes and superconvergence of their associated long-time error has also been investigated by \cite{asthana2017consistency}. The error estimate derived by \cite{asthana2017consistency} is however, semi-analytical and requires the numerical evaluation of the order of accuracy of the dispersion-dissipation error. Although an error bound for the L2-error of the ESFR spatial discretization has been formally derived by \cite{lambert2023l2}, to the knowledge of the authors of this paper, there is no existing work focusing on proving formally the superconvergent properties of ESFR schemes. The objective of this paper is thus to address this theoretical gap. Namely, we present a simple proof for the superconvergence of the combined dispersion-dissipation error of ESFR schemes for the linear advection problem when using an upwind numerical flux. Our results provide firm theoretical grounds to explain the numerical observations of \cite{vincent2011insights, vermeire2016properties} and when considered in conjunction with the work of \cite{asthana2017consistency}, they allow for a fully analytic characterization of the superconvergence of the long-time error of ESFR schemes.

The paper is organized as follows. \Cref{sec:preliminaries} provides a brief review of the FR framework and its relationship with ESFR schemes. \Cref{sec:proof} focuses on deriving a theoretical dispersion-dissipation error estimate for ESFR schemes and hence proving the superconvergent properties of the latter. In \cref{sec:numexp}, these theoretical results are then validated via numerical experiments. Finally, in \cref{sec:extension}, extension of the theoretical results to more general symmetric FR schemes is briefly considered.

\section{Preliminaries}
\label{sec:preliminaries}

\subsection{ESFR schemes for the Linear Advection Problem}
In this section, we briefly review the ESFR framework in the context of the one-dimensional linear advection problem. For a more thorough presentation of the matter, we refer the reader to the work of \cite{huynh2007flux,vincent2011new,allaneau2011connections}.

For a unit propagation speed, the linear advection problem is written as
\begin{equation}
    \frac{\partial u}{\partial t} + \frac{\partial u}{\partial x} = 0
    \quad ; \quad
    u(x,0) = u_0(x)
    \quad ; \quad
    x \in \Omega,
    \label{eq:conservation_law}
\end{equation}
where $u$ is the conserved scalar and $\Omega$ is a closed interval in $\mathbb{R}$. Following the classical FR formulation \cite{huynh2007flux}, we start by discretizing $\Omega$ into $N$ non-overlapping elements $\Omega_k$ such that
\begin{equation}
    \Omega = \bigcup_{k=0}^{N-1} \Omega_k
    \quad ; \quad
    \bigcap_{k=0}^{N-1} \Omega_k = \emptyset.
    \label{eq:discretization}
\end{equation}
Upon introducing the standard normalization
\begin{equation}
    \hat{u}_k(\xi,t) := u(\Gamma_k^{-1}(\xi), t)
\end{equation}
via the elemental linear mappings $\Gamma_k : \Omega_k \to [-1,1]$
\begin{equation}
    \xi = \Gamma_k(x) = 2\frac{x-x_k}{x_{k+1}-x_k} -1,
\end{equation}
and their associated inverses
\begin{equation}
    x=\Gamma^{-1}_k(\xi) = \frac{1-\xi}{2}x_k + \frac{1+\xi}{2}x_{k+1},
\end{equation}
we can rewrite \cref{eq:conservation_law} as
\begin{equation}
    \frac{\partial \hat{u}_k}{\partial t} + \frac{2}{|\Omega_k|}\frac{\partial \hat{u}_k}{\partial \xi} = 0.
    \label{eq:conservation_law_ref}
\end{equation}
In the FR framework, the local solutions $\hat{u}_k$ are then approximated using a combination of $p+1$ linearly independent high-order polynomials $\chi_j \in \mathcal{P}^{p}([-1,1])$ such that
\begin{equation}
    \hat{u}_k(\xi,t) \approx \hat{u}_k^h(\xi,t) = \sum_{j=0}^{p} \hat{u}^h_{jk}(t) \chi_j(\xi) \label{eq:local_approx}
\end{equation}
and the global solution on $\Omega$ is consequently given by
\begin{equation}
    u(x,t) \approx u^h(x,t) = \bigoplus_{k=0}^{N-1} \hat{u}_k^h(\Gamma_k(x),t).
\end{equation}
The classical FR spatial discretization is then obtained by requiring the local solutions introduced in \cref{eq:local_approx} to satisfy \cref{eq:conservation_law_ref} in its differential form and by introducing the left and right correction functions $g_L, g_R \in \mathcal{P}^{p+1}([-1,1])$ to ensure continuity of the flux at the left and right elemental interfaces. This yields the familiar semi-discrete scheme
\begin{equation}
    \sum_{j=0}^{p} \frac{d\hat{u}_{jk}^h}{dt}\chi_j
    + \frac{2}{|\Omega_k|}\left(\sum_{j=0}^{p}\hat{u}_{jk}^h \frac{d\chi_j}{d\xi}
    + (\hat{f}^*_{kL} - \hat{u}^h_k|_{L})\frac{dg_L}{d\xi} + (\hat{f}^*_{kR} - \hat{u}^h_k|_{R})\frac{dg_R}{d\xi}\right) = 0,
    \label{eq:FR}
\end{equation}
where $f^*_{kL}$ and $f^*_{kL}$ denote the numerical flux at the left and right elemental interfaces and $\hat{u}^h_k|_{L} = \hat{u}^h_k(-1,t)$ and $\hat{u}^h_k|_{R} = \hat{u}^h_k(1,t)$. In vector form, \cref{eq:FR} amounts to
\begin{equation}
    \frac{d \hat{\mathbf{u}}^h_k}{dt}
    +
    \mathbf{D} \hat{\mathbf{f}}^D_k
    + (\hat{f}^*_{kL} - \mathbf{l}^{T} \hat{\mathbf{f}}^D_k) \mathbf{h}_L
    + (\hat{f}^*_{kR} - \mathbf{r}^{T} \hat{\mathbf{f}}^D_k) \mathbf{h}_R =0,
    \label{eq:FR_operators}
\end{equation}
where
\begin{equation}
    (\mathbf{l})_j := \chi_j(-1)
    \quad ; \quad
    (\mathbf{r})_j = \chi_j(1)
    \quad ; \quad
    (\mathbf{u}_k^h)_j = \hat{u}^h_{jk},
\end{equation}
\begin{equation}
    g_L'(\xi) = \sum_{j=0}^p (\mathbf{h}_L)_j \chi_j(\xi)
    \quad ; \quad
    g_R'(\xi) = \sum_{j=0}^p (\mathbf{h}_R)_j \chi_j(\xi)
\end{equation}
and $\mathbf{D}$ is the matrix representation of the differential operator on $\mathcal{P}^{p}([-1,1])$.

ESFR schemes were initially derived by \cite{vincent2011new} by requiring $g_L$ and $g_R$ to be such that $u^h$ is stable with respect to a broken Sobolev norm. \cite{allaneau2011connections,zwanenburg2016equivalence} showed that this was equivalent to constructing $\mathbf{h}_L$ and $\mathbf{h}_R$ via
\begin{equation}
    \mathbf{h}_R = (\mathbf{M} + \mathbf{K})^{-1} \mathbf{r}
    \quad ; \quad
    \mathbf{h}_L = -(\mathbf{M} + \mathbf{K})^{-1}\mathbf{l},
\end{equation}
where $\mathbf{M}$ is the mass matrix and $\mathbf{K}$ is the ESFR filtering matrix defined as
\begin{equation}
    \mathbf{M}_{mn} := \int_{-1}^{1} \chi_m \chi_n d\xi
    \quad ; \quad
    \mathbf{K} := \frac{1}{2}c(\mathbf{D}^{p})^T\mathbf{M}\mathbf{D}^p
\end{equation}
and $c$ is the ESFR scalar parameter. Clearly, when $c=0$, a standard nodal discontinuous Galerkin scheme is recovered. As shown by \cite{vincent2011new}, ESFR schemes are linearly stable provided that $c > c_{-}$. Moreover, when $c = c_{SD}$ or $c = c_{HU}$, one recovers a spectral difference scheme or Huynh's $g_2$ scheme, respectively. We refer the reader to the work of \cite{vincent2011new} for the specific analytical expressions for $c_{-}$, $c_{SD}$ and $c_{HU}$ and for a more thorough description of ESFR schemes.

\subsection{Dispersion-dissipation Error of ESFR Schemes}
As previously foreshadowed, this study focuses on characterizing the superconvergent properties of the dispersion-dissipation error of ESFR schemes. In this context, consistent with the work of \cite{hesthaven2007nodal,vincent2011insights,vermeire2016properties}, dispersion-dissipation error is defined via
\begin{equation}
    E_T(\theta) := |\omega^h(\theta) - \omega_{PDE}(\theta)|,
    \label{eq:spectral_error}
\end{equation}
where $\omega^h: \mathbb{R} \to \mathbb{C}$ is the dispersion-dissipation relation of the numerical scheme, $\omega_{PDE}: \mathbb{R} \to \mathbb{C}$ is the dispersion-dissipation relation of the PDE solved by the numerical scheme and $\theta \in \mathbb{R}$ represents the non-dimensional wavenumber of Bloch wave solutions. For the linear advection problem (unit propagation speed), the true dispersion-dissipation relation is given by
\begin{equation}
    \omega_{PDE}(\theta) = \theta,
\end{equation}
which ensures Bloch wave solutions of the form
\begin{equation}
    u(x,t) = \exp(i(\theta x - \omega_{PDE}(\theta)t))
\end{equation}
are solutions to \cref{eq:conservation_law} and where $i := \sqrt{-1}$. Following the same approach as \cite{vincent2011insights,moura2015linear,vermeire2016properties}, the dispersion-dissipation relations of ESFR schemes are in turn obtained by substituting numerical Bloch waves of the form
\begin{equation}
    \hat{\mathbf{u}}^h_k = \exp (i(k\theta - \omega^h(\theta)t)) \hat{\mathbf{v}}
    \label{eq:Bloch_wave}
\end{equation}
in \cref{eq:FR_operators}, which yields the eigenvalue problem
\begin{equation}
    \mathbf{H}(\theta) \hat{\mathbf{v}}(\theta) = \omega^h(\theta) \hat{\mathbf{v}}(\theta).
    \label{eq:Bloch_eig}
\end{equation}
In what follows, we will refer to $\mathbf{H}$ as the von Neumann matrix of a given numerical scheme. For an upwind numerical flux, \cite{vincent2011insights} showed that $\mathbf{H} : \mathbb{R} \to \mathbb{C}^{(p+1)\times(p+1)}$ is given by
\begin{equation}
    \mathbf{H}(\theta) = -2i\left(
    \mathbf{D} + (\mathbf{M} + \mathbf{K})^{-1} \mathbf{l}\mathbf{l}^T \nonumber - (\mathbf{M} + \mathbf{K})^{-1} \mathbf{l} \mathbf{r}^T  e^{-i\theta}\right),
\end{equation}
and $\omega^h(\theta)$ is defined as the unique physical eigenvalue of $\mathbf{H}(\theta)$ which approximates the true dispersion-dissipation relation of the PDE as $\theta \to 0$.

As discussed by \cite{huynh2007flux}, studying the convergence of $E_T$ as $\theta \to 0$ allows for the characterization of the impact of $h$-refinement on the dispersion-dissipation error of a numerical scheme. Following the notation of \cite{huynh2007flux,vincent2011insights,vermeire2016properties}, we define the spectral order of accuracy of a numerical scheme via
\begin{equation}
    E_T \propto \theta^{A_T + 1} \quad \text{as} \quad \theta \to 0.
\end{equation}
In other words, $A_T$ is the convergence rate of $E_T$ minus one as $\theta \to 0$. ESFR schemes are said to exhibit superconvergent behavior because $A_T > p$ for this family of linearly stable schemes. The superconvergent properties of FR schemes have been observed by \cite{huynh2007flux,vincent2011insights,vermeire2016properties} through numerical experiments. More precisely, by selecting a sufficiently small value $0 < \Delta \theta  \ll 1$, they computed $E_T$ numerically and used
\begin{equation}
    A_T \approx \frac{\ln(E_T(\Delta \theta)) - \ln(E_T(\Delta \theta / 2))}{\ln(2)} - 1
    \label{eq:spectral_OA}
\end{equation}
to estimate the convergence rate of the dispersion-dissipation error. They found that $A_T \approx 2p+1$ for values of $c$ close to 0, that $A_T \approx 2p$ for $c \in c_{SD}, c_{HU}$ and that $A_T \approx 2p-1$ for sufficiently large values of $c$. In the rest of this study, through analytical means, we will determine $A_T$ rigorously for ESFR schemes and use this to unveil the fundamental mechanism behind superconvergence of FR methods.

\section{Proof for Superconvergence Property}
\label{sec:proof}

\subsection{Equivalent Problem}
As noted by \cite{asthana2017consistency, huynh2007flux}, the challenge of showing the superconvergent behavior of ESFR schemes essentially lies in the difficulty of finding an analytical expression for the physical eigenvalue of $\mathbf{H}$ when $p > 2$. In this work, to circumvent this issue, we first introduce \cref{theo:equivalence}, which connects the order of convergence of $E_T$ to that of the magnitude of the determinant of a modified von Neumann matrix. We will then proceed to showing the superconvergent properties of this determinant.
\begin{theorem}
\label{theo:equivalence}
Let $\tilde{\mathbf{H}} : \mathbb{R} \to \mathbb{C}^{(p+1)\times(p+1)}$ be defined as
\begin{equation}
    \tilde{\mathbf{H}}(\theta) :=(\mathbf{M} + \mathbf{K})\left(\mathbf{D} -\frac{1}{2}i\theta \mathbf{I}\right) +\mathbf{l}\left(\mathbf{l}-\mathbf{r}e^{-i\theta}\right)^T.
    \label{eq:tildeH}
\end{equation}
Then, as $\theta \to 0$, $E_T(\theta)$ and $|\text{\textup{det}}(\tilde{\mathbf{H}}(\theta))|$ converge to 0 at the same rate provided that the ESFR parameter satisfies $c_{-} < c < \infty$.
\end{theorem}
\begin{proof}
    First, note that by definition $\tilde{\mathbf{H}}(\theta)$ is related to $\mathbf{H}(\theta)$ via
    \begin{equation}
        \tilde{\mathbf{H}}(\theta) = -\frac{1}{2i}(\mathbf{M}+\mathbf{K})(\mathbf{H}(\theta) - \theta\mathbf{I}).
    \end{equation}
    Hence, we must have
    \begin{equation}
        \text{det}(\tilde{\mathbf{H}}(\theta)) =
        \left(-\frac{1}{2i}\right)^{p+1}\text{det}(\mathbf{M}+\mathbf{K}) \text{det}(\mathbf{H}(\theta) - \theta \mathbf{I})),
        \label{eq:modified_VonNeumann}
    \end{equation}
    where $\text{det}(\mathbf{M}+\mathbf{K}) \neq 0$ provided that $c$ satisfies the prescribed conditions. Moreover, it can readily be seen that
    \begin{align}
        |\text{det}(\mathbf{H}(\theta) - \theta\mathbf{I})|
        &=
        |\omega^h(\theta) - \theta||\lambda_1(\theta) - \theta| \cdots |\lambda_p(\theta) - \theta| \\
        &= 
        E_T |\lambda_1(\theta) - \theta| \cdots |\lambda_p(\theta) - \theta|,
    \end{align}
    where $\omega^h(\theta), \lambda_1(\theta),...,\lambda_p(\theta)$ are the $p+1$ eigenvalues of $\mathbf{H}(\theta)$. As shown in \cref{app:lem_mult}, the algebraic multiplicity of the eigenvalue $\omega^h(0) = 0$ is one provided that $c_{-} < c < \infty$. Additionally, since the entries of $\mathbf{H}$ depend continuously on $\theta$, each eigenvalue of $\mathbf{H}$ is a continuous function of $\theta$ \cite{horn2013roger} and $\omega^h$ is the unique eigenvalue approaching $0$ as $\theta \to 0$. Thus, $\exists C_1, C_2 > 0$ such that
\begin{equation}
    C_1 |\text{det}(\tilde{\mathbf{H}}(\theta))| \leq E_T(\theta) \leq C_2|\text{det}(\tilde{\mathbf{H}}(\theta))|
\end{equation}
for all values of $\theta$ sufficiently small.
\end{proof}

\begin{remark}
    One should note that since the algebraic multiplicity of $\omega^h$ is one in the vicinity of $\theta=0$, by the implicit function theorem, $\omega^h$ is $C^{\infty}$ at $\theta = 0$ and thus the order of convergence of $E_T$ can be safely defined.
\end{remark}

\subsection{Superconvergence of the Dispersion-Dissipation Error}
Based on \cref{theo:equivalence}, we now proceed to deriving an analytical expression for $|\text{det}(\tilde{\mathbf{H}}(\theta))|$. We begin by showing that the convergence behavior of the latter is related to a specific rational approximation of the exponential function, as per \cref{lem:rational_approx}.
\begin{lemma}
\label{lem:rational_approx}
Let $\psi^{(k)}_j(\xi)$ denote the $k$th derivative of the $j$th Legendre polynomial evaluated at $\xi$ and let $q_j$ denote the $j$th diagonal entry of $(\mathbf{M}+\mathbf{K})^{-1}$ when expressed in the Legendre basis. The determinant of the modified von Neumann matrix $\tilde{\mathbf{H}}$ is given by
\begin{equation}
     \text{\textup{det}}(\tilde{\mathbf{H}}(\theta)) =
    (-1)^{p+1}\text{\textup{det}}(\mathbf{M}+\mathbf{K})
    \left(P(i\theta) - Q(i\theta)e^{-i\theta}\right),
    \label{eq:simp1}
\end{equation}
where $P \in \mathcal{P}^{p+1}$ and $Q \in \mathcal{P}^{p}$ are defined as
\begin{align}
    P(z) &:=2^{-(p+1)}z^{p+1}-\sum_{j,k=0}^p (-1)^j2^{k-p} q_j\psi^{(k)}_j(-1)z^{p-k}, \label{eq:approximant_P} \\
    Q(z) &:=-\sum_{j,k=0}^p (-1)^j2^{k-p} q_j\psi^{(k)}_j(1)z^{p-k}. \label{eq:approximant_Q}
\end{align}
\end{lemma}
\begin{proof}
    For convenience and without loss of generality, we will assume that $\tilde{\mathbf{H}}$ is constructed in the Legendre basis. Leveraging the dyadic product structure of $\tilde{\mathbf{H}}$, we start by applying the well-known matrix determinant lemma to \cref{eq:tildeH}, which yields
\begin{equation}
    \text{det}(\tilde{\mathbf{H}}(\theta)) =
    \text{det}(\mathbf{A}(\theta))
    \left(1 + \mathbf{l}^T\mathbf{A}(\theta)^{-1} \mathbf{l} - \mathbf{r}^T\mathbf{A}(\theta)^{-1}\mathbf{l}e^{-i\theta}\right),
    \label{eq:det_lemma}
\end{equation}
where
\begin{equation}
    \mathbf{A}(\theta) := (\mathbf{M} + \mathbf{K})\left(\mathbf{D} -\frac{1}{2}i\theta \mathbf{I}\right).
\end{equation}
Since $\mathbf{D}$ is strictly upper triangular, we can deduce
\begin{equation}
    \text{det}(\mathbf{A}(\theta))
    = 
    \text{det}(\mathbf{M} +\mathbf{K})\left(-\frac{1}{2}i\theta\right)^{p+1}.
\end{equation}
To determine $\mathbf{A}(\theta)^{-1}$, we expand it in terms of its von Neumann series. In particular, since $\mathbf{D}$ is nilpotent with $\mathbf{D}^{p+1} = 0$, the von Neumann series of this operator can be written as the finite sum
\begin{align}
    \mathbf{A}(\theta)^{-1}
    &= -\left(\frac{1}{2}i\theta \mathbf{I} - \mathbf{D}\right)^{-1}(\mathbf{M}+\mathbf{K})^{-1} \\
    &=
    -\sum_{k=0}^{p} 2^{k+1}(i\theta)^{-(k+1)} \mathbf{D}^k (\mathbf{M}+\mathbf{K})^{-1}.
    \label{eq:vonneumann_series}
\end{align}
Combining \cref{eq:det_lemma} and \cref{eq:vonneumann_series} yields
\begin{equation}
     \text{det}(\tilde{\mathbf{H}}(\theta)) =
    (-1)^{p+1}\text{det}(\mathbf{M}+\mathbf{K})
    \left(2^{-(p+1)}(i\theta)^{p+1} + \mathbf{l}^T\mathbf{B}(\theta) \mathbf{l} - \mathbf{r}^T\mathbf{B}(\theta)\mathbf{l}e^{-i\theta}\right),
\end{equation}
where
\begin{equation}
    \mathbf{B}(\theta) := \sum_{k=0}^{p} 2^{k-p}(i\theta)^{p-k} \mathbf{D}^k (\mathbf{M}+\mathbf{K})^{-1}.
\end{equation}
Letting $q_j$ denote the $j$th diagonal entry of $(\mathbf{M}+\mathbf{K})^{-1}$ we leverage the structure of $\mathbf{B}$ and the fact $\mathbf{M} + \mathbf{K}$ is diagonal to find that
\begin{align}
    \mathbf{l}^T\mathbf{B}(\theta)\mathbf{l}
    &=
    -\sum_{k=0}^{p} 2^{k-p}(i\theta)^{p-k} \mathbf{l}^T\mathbf{D}^k (\mathbf{M}+\mathbf{K})^{-1} \mathbf{l} \\
    &=
    -\sum_{j=0}^p\sum_{k=0}^{p} (-1)^j2^{k-p} q_j\psi^{(k)}_j(-1)(i\theta)^{p-k}.
    \label{eq:left}
\end{align}
Similarly, it can be seen that
\begin{align}
    \mathbf{r}^T\mathbf{B}(\theta)\mathbf{l}
    &=
    -\sum_{k=0}^{p} 2^{k-p}(i\theta)^{p-k} \mathbf{r}^T\mathbf{D}^k (\mathbf{M}+\mathbf{K})^{-1} \mathbf{l} \\
    &=
    -\sum_{j=0}^p\sum_{k=0}^{p} (-1)^j2^{k-p} q_j\psi^{(k)}_j(1)(i\theta)^{p-k}.
    \label{eq:right}
\end{align}
\end{proof}
From \cref{lem:rational_approx}, one can see that the determinant of the modified von Neumann matrix $\tilde{\mathbf{H}}$ will exhibit superconvergence properties as $\theta \to 0$ if and only if $P(z) / Q(z)$ is a superconvergent rational approximant of $e^{-z}$ about $z=0$. To prove the latter, we consider generalizing the work of \cite{ahmad1998orthogonal}. Namely, \cite{ahmad1998orthogonal} has shown a simple connection between Legendre polynomials and $(p,p)$ Padé approximants of the exponential function. Here, we extend this approach to $(p+1,p)$ Padé approximants and use this to show that the polynomials defined in \cref{eq:approximant_P} and \cref{eq:approximant_Q} form a general class of superconvergent rational approximants of $e^{-z}$. Our key findings are summarized in \cref{lem:Padé}.
\begin{lemma}
\label{lem:Padé}
The rational function $P(z)/Q(z)$ with $P$ and $Q$ given by \cref{eq:approximant_P} and \cref{eq:approximant_Q} is a superconvergent approximant of $e^{-z}$. Namely,
\begin{align}
    P(z) - Q(z)e^{-z}
    &=
    \frac{1}{2^{p+1}}z^{p+1}e^{-z/2}\left(e^{z/2} - \sum_{j=0}^p2q_jI_j(z/2)\right) \label{eq:Padé1} \\
    &= \frac{1}{2^{p+1}}
    \left(b_1(p,c)z^{2p+1} + b_2(p,c)z^{2p+2}\right) + O(z^{2p+3}) \label{eq:Padé2}
\end{align}
where
\begin{align}
    b_1(p,c) &:= \frac{(1-f(c))p!}{(2p)!}, \label{eq:b1} \\
    b_2(p,c) &:= \frac{(f(c)-1)p!}{2(2p)!} + \frac{(p+1)!}{(2p+2)!}, \label{eq:b2} \\
    f(c) &:= \frac{1}{1+c\frac{(2p+1)(2p!)^2}{2^{2p+1}(p!)^2}}
\end{align}
and $I_j$ stands for the $j$th modified spherical Bessel function of the first kind.
\end{lemma}
\begin{proof}
We start with the well-known series expansion of the exponential function in the Legendre basis \cite{NIST:DLMF},
\begin{equation}
    e^{y \cos(\alpha)}
    =
    \sum_{j=0}^{\infty}(2j+1)I_j(y) \psi_j(\cos(\alpha)),
    \label{eq:master_identity}
\end{equation}
where $y \in \mathbb{C}$ and $\alpha \in \mathbb{R}$. Setting $\alpha = 0$ in \cref{eq:master_identity} yields the series expansion
\begin{equation}
    e^y = \sum_{j=0}^{\infty}(2j+1)I_j(y).
    \label{eq:identity2}
\end{equation}
Moreover, letting $x:=\cos(\alpha)$, \cref{eq:master_identity} can be used to compute the projection of the exponential function on $\psi_j$
\begin{equation}
    \int_{-1}^1 e^{yx} \psi_j(x)dx
    =
    \sum_{k=0}^{\infty}(2k+1)I_j(y) \int_{-1}^1 \psi_j(x)\psi_k(x)dx
    =
    2I_j(y),
    \label{eq:identity1}
\end{equation}
from which one can see that
\begin{equation}
    \int_{-1}^1 e^{-yx} \psi_j(x)dx
    =
    2I_j(-y)
    =2(-1)^jI_j(y).
    \label{eq:aux1}
\end{equation}
Following the approach of \cite{ahmad1998orthogonal}, we now consider integrating by parts the left-hand side of \cref{eq:aux1} $j+1$ times, which yields
\begin{equation}
    \left[-e^{-yx}\sum_{k=0}^j \frac{\psi_j^{(k)}(x)}{y^{k+1}}\right]_{-1}^{1}
    =
    \left[-e^{-yx}\sum_{k=0}^p \frac{\psi_j^{(k)}(x)}{y^{k+1}}\right]_{-1}^{1}
    =
    2(-1)^jI_j(y).
    \label{eq:aux2}
\end{equation}
Upon multiplication of both sides of \cref{eq:aux2} by $(-1)^{j+1}y^{p+1}q_je^{-y}$, one obtains
\begin{equation}
    -\sum_{k=0}^p (-1)^j q_j \psi_j^{(k)}(-1) y^{p-k}
    +
    e^{-2y}\sum_{k=0}^p (-1)^j q_j \psi_j^{(k)}(1) y^{p-k}
    =
    -2 q_j e^{-y}y^{p+1}I_j(y),
    \label{eq:aux3}
\end{equation}
which we then sum from $j=0$ to $j=p$ to get
\begin{align}
    &-\sum_{j=0}^{p}\sum_{k=0}^p(-1)^jq_j \psi_j^{(k)}(-1) y^{p-k} \nonumber
    +e^{-2y}\sum_{j=0}^p\sum_{k=0}^p (-1)^jq_j\psi_j^{(k)}(1) y^{p-k} \\
    &=
    -2e^{-z}y^{p+1}\sum_{j=0}^pq_jI_j(y).
    \label{eq:aux4}
\end{align}
Letting $y = z/2$, adding $2^{-(p+1)}z^{p+1}$ and comparing with \cref{eq:approximant_P} and \cref{eq:approximant_Q}, \cref{eq:aux4} becomes
\begin{align}
    &\frac{1}{2^{p+1}}z^{p+1}-\sum_{j,k=0}^p(-1)^j 2^{k-p} q_j \psi_j^{(k)}(-1) z^{p-k} \nonumber
    +e^{-z}\sum_{j,k=0}^p (-1)^j 2^{k-p} q_j\psi_j^{(k)}(1) z^{p-k} \\
    &=
    P(z) - Q(z)e^{-z} \nonumber \\
    &=
    \frac{1}{2^{p+1}}z^{p+1}e^{-z/2}\left(e^{z/2} - \sum_{j=0}^p2q_jI_j(z/2)\right),
    \label{eq:aux5}
\end{align}
which completes the proof for \cref{eq:Padé1}. From \cref{eq:aux5}, one can see that the rational approximant $P(z)/Q(z)$ is at least $O(z^{p+1})$ accurate. To show superconvergence of the approximant, the rightmost factor in \cref{eq:aux5} must, however, be studied more closely.

For ESFR schemes, it has been shown that the diagonal entries of $(\mathbf{M} + \mathbf{K})^{-1}$ were given by \cite{allaneau2011connections, zwanenburg2016equivalence}
\begin{equation}
q_j =
    \begin{cases}
        \frac{2j+1}{2}& \quad \text{for} \quad 0 \leq j< p, \\
        \frac{2j+1}{2}f(c)& \quad \text{for} \quad j=p,
    \end{cases}
\end{equation}
where
\begin{equation}
    f(c) := \frac{1}{1+c\frac{(2p+1)(2p!)^2}{2^{2p+1}(p!)^2}}.
\end{equation}
Hence, it follows that when $c=0$ -- which corresponds to a DG scheme -- the rightmost factor in \cref{eq:aux5} is exactly equal to the error incurred when the series expansion given by \cref{eq:identity2} is truncated at the $p$th order. Since the associated truncation error is $O(z^{p+1})$ in this case, the rational approximant $P(z)/Q(z)$ is superconvergent with order $O(z^{2p+2})$. When $c \neq 0$, one additional term in the partial sum is lost, and the truncation error is now $O(z^{p})$, yielding a superconvergent rational approximant with order $O(z^{2p+1})$.

To obtain an explicit expression for the leading term of the error associated with the rational approximant $P(z)/Q(z)$, we first recall that the $j$th spherical Bessel function of the first kind is defined via the power series
\begin{equation}
    I_j(z) = \frac{\sqrt{\pi}}{2} \sum_{k=0}^{\infty} \frac{\left(\frac{1}{2}z\right)^{2k+j}}{k! \Gamma\left(k+j+\frac{3}{2}\right)},
    \label{eq:Bessel}
\end{equation}
where $\Gamma$ represents the Gamma function. Using \cref{eq:aux5}, \cref{eq:identity2} and \cref{eq:Bessel}, one can see that
\begin{align}
    &P(z) - Q(z)e^{-z} =\frac{1}{2^{p+1}}z^{p+1}e^{-z/2}\left(e^{z/2} - \sum_{j=0}^p2q_jI_j(z/2)\right) \nonumber \\
    &=
    (1-f(c))a_{p}(z/2)^{2p+1}
    + ((f(c)-1)a_{p}+a_{p+1})(z/2)^{2p+2} + O(z^{2p+3}),
\end{align}
where
\begin{equation}
    a_{p} := \frac{\sqrt{\pi}}{2} \frac{2p+1}{2^p \Gamma\left(p-\frac{3}{2}\right)}.
\end{equation}
Using the properties of the Gamma function, this can be simplified to
\begin{equation}
    a_{p} = \frac{2^pp!}{(2p)!}.
\end{equation}
Combining these results, we finally get
\begin{equation}
    P(z)-Q(z)e^{-z} =
    \left(b_1(p,c)z^{2p+1} + b_2(p,c)z^{2p+2}\right) + O(z^{2p+3}),
\end{equation}
which completes the proof for \cref{eq:Padé2}.
\end{proof}
As can be seen from \cref{lem:Padé}, the DG scheme is of particular interest when considering the class of superconvergent rational approximants associated with ESFR schemes. Namely, when $c=0$, $f(c)=1$ and $b_1(p,c) = 0$, yielding a rational approximant which is truly $O(z^{2p+2})$. In this case, $P(z)/Q(z)$ exactly coincides with the $(p+1,p)$ Padé approximant of the function $e^{-z}$. This can readily be shown through careful examination of the coefficients of $P(z)$ and $Q(z)$ or by recalling that the $(m,n)$ Padé approximation of a function is the unique rational approximant with a truncation error of order $O(z^{m+n+1})$. As expected, this Padé approximant is identical to that used in the work of \cite{ainsworth2004dispersive,atkins2009super,chalmers2017spatial} to show superconvergent properties of DG schemes.

Combining the results of \cref{lem:rational_approx}, \cref{lem:Padé} and \cref{theo:equivalence}, we are now ready to state the main theoretical result of this paper, which is shown in \cref{theo:error_estim}.
\begin{theorem}
\label{theo:error_estim}
Consider an arbitrary ESFR scheme with $c_{-} < c < \infty$. Then, $\exists C_1, C_2 > 0$ (dependent on $c$) such that for all $\theta$ sufficiently small
\begin{equation}
    C_1F(p,c,\theta) \leq E_T(\theta) \leq C_2F(p,c,\theta),
\end{equation}
where
\begin{equation}
    F(p,c,\theta) :=|b_1(p,c)(i\theta)^{2p+1} + b_2(p,c)(i\theta)^{2p+2}|
\end{equation}
and $b_1$ and $b_2$ are given by \cref{eq:b1} and \cref{eq:b2} respectively. 
\end{theorem}
\begin{proof}
    \Cref{theo:error_estim} follows from direct application of \cref{lem:rational_approx}, \cref{lem:Padé} and \cref{theo:equivalence}.
\end{proof}

\subsection{Important Remarks}
In the previous section, we essentially showed that as $\theta \to 0$, the dispersion-dissipation error of ESFR schemes scaled as per
\begin{equation}
    E_T(\theta) \propto |b_1(p,c)(i\theta)^{2p+1} + b_2(p,c)(i\theta)^{2p+2}|.
    \label{eq:error_estim}
\end{equation}
There is a seeming contradiction between this theoretical error estimate and the experimental results of \cite{vincent2011insights}. Indeed, while \cref{eq:error_estim} does predict that $A_T$ should drop from $A_T\approx2p+1$ to $A_T\approx2p$ as $c$ is progressively increased from $c=0$, it does not explain why a subsequent drop from $A_T\approx2p$ to $A_T\approx2p-1$ is observed at even larger values of $c$.

This apparent contradiction can be resolved by noting that our theoretical error estimate is not valid in the limit as $c \to \infty$. Indeed, as shown in \cref{app:lem_mult}, when $c \to \infty$, $\mathbf{M} + \mathbf{K}$ becomes singular and as a result, $\text{det}(\tilde{\mathbf{H}})$ vanishes. Thus, although the eigenvalue $\omega^h(0)=0$ has multiplicity one for any finite value of $c$, in the limit, its multiplicity is actually two, which invalidates the assumptions of \cref{theo:equivalence}. To resolve this issue, we can simply remove the problematic eigenvalue by dividing $\text{det}(\tilde{\mathbf{H}})$ by the latter. Namely, we write
\begin{align}
    E_T
    &\propto
    |\omega^h(\theta,c)-\theta||\lambda_2(\theta,c) - \theta|\cdots |\lambda_p(\theta,c)-\theta|
    \propto
    \frac{|\text{det}(\tilde{\mathbf{H}}(\theta,c))|}{|\lambda_1(\theta,c)-\theta|} \nonumber \\
    &\propto
    \frac{\left|b_1(p,c)(i\theta)^{2p+1} + b_2(p,c)(i\theta)^{2p+2}\right|}{|\lambda_1(\theta,c)-\theta|},
    \label{eq:final}
\end{align}
where $\lambda_1$ is the unique eigenvalue of $\mathbf{H}$ with the properties
\begin{equation}
    \forall c \in \mathbb{R} : \lambda_1(0,c) \neq 0
    \quad \text{and} \quad
    \lim_{c \to \infty} \lambda_1(0,c) = 0.
    \label{eq:eig_special}
\end{equation}
Note that in \cref{eq:final} and \cref{eq:eig_special}, the dependence of the eigenvalues of $\mathbf{H}$ on $\theta$ and $c$ has been explicitly written out for clarity. Careful observation of \cref{eq:final} leads to interesting conclusions. First, as can be seen, unlike \cref{eq:error_estim}, \cref{eq:final} is valid in the limit as $c \to \infty$. Namely, upon computing the limit, one finds that $E_T$ is of order $O(\theta^{2p})$ and thus that $A_T = 2p-1$, which agrees with the fact that ESFR schemes collapse to a DG scheme of degree $p-1$ when $c \to \infty$. It can also be observed that for any finite value of $c$, since the denominator of \cref{eq:final} does not vanish, \cref{eq:error_estim} and \cref{eq:final} are equivalent and the orders of accuracy predicted by \cref{eq:error_estim} should be regarded as correct.

\Cref{eq:final}, however, does provide a simple explanation for the discrepancy observed between our theoretical error estimate and the results of \cite{vincent2011insights}. Namely, if one attempts to estimate $A_T$ numerically via \cref{eq:spectral_OA} at large values of $c$, it is clear from \cref{eq:final} that the small magnitude of $\lambda_1(c,0)$ will make it very difficult to measure experimentally the true asymptotic convergence rate. In particular, if $A_T$ is not computed at a sufficiently small value of $\Delta\theta$, the denominator of \cref{eq:final} will behave as if it were $O(\theta)$ and the convergence rate will lose one additional order. Given that the values of $\Delta\theta$ selected by \cite{vincent2011insights} were $\pi/8$, $\pi/4$, $\pi/3$ and $2\pi/5$ when $p$ was 2,3,4 and 5, respectively, it appears likely that the convergence rates they measured were not in the asymptotic regime. This was confirmed numerically as shown in the next section.

\section{Numerical Experiments}
\label{sec:numexp}
This section presents the numerical experiments performed to confirm the validity of the dispersion-dissipation error estimate introduced in the previous section.

\subsection{Convergence of the Dispersion-Dissipation Error}
\begin{figure}
    \centering
    \includegraphics[width=0.6\textwidth]{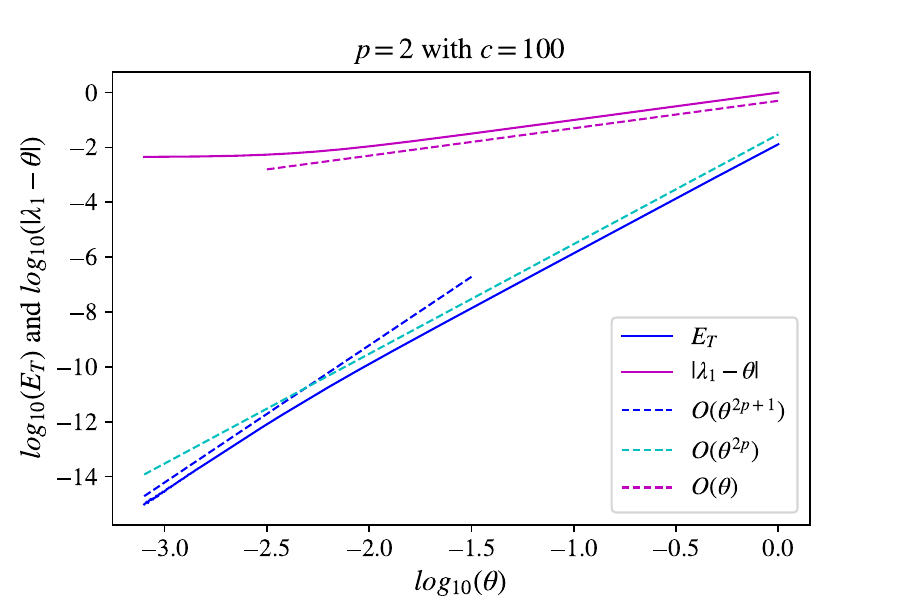}
    \caption{Convergence rate of $E_T$ and $|\lambda_1(c,\theta)-\theta|$ for $p=2$ and $c=100$.}
    \label{fig:myproof}
\end{figure}

In the previous section, it was hypothesized that the values of $A_T$ reported by \cite{vincent2011insights} for large values of $c$ were actually not measured in the asymptotic regime. To verify this claim, $E_T$ was computed numerically via \cref{eq:spectral_error} and \cref{eq:Bloch_eig} and plotted against $\theta$ for the ESFR scheme characterized by $p=2$ and $c=100$. For this scheme, \cite{vincent2011insights} reports that $A_T = 2p-1$, which implies $E_T \propto \theta^{2p}$. Results are shown in \cref{fig:myproof}. As can be seen, it is clear that although it initially appears that $E_T \propto \theta^{2p}$, if $\theta$ is taken to be sufficiently small, the theoretical convergence rate of $E_T \propto \theta^{2p+1}$ can be observed.

\Cref{fig:myproof} also shows the variation of the denominator of \cref{eq:final} with respect to $\theta$. As can be seen, for sufficiently small values of $\theta$, $\lambda_1$ approaches a non-zero plateau. Hence, the algebraic multiplicity of the eigenvalue $\omega^h$ is one in the vicinity of $\theta = 0$, as expected. Consistent with \cref{eq:final}, the change in the rate of decay of $E_T$ coincides with the stabilization in the value of $|\lambda_1(\theta,c) - \theta|$.

It should finally be noted that this numerical experiment becomes prohibitively difficult to conduct when the scheme order is large. Indeed, as a result of the superconvergent behavior of $E_T$, for large values of $p$, standard double machine precision is quickly reached before $E_T$ enters the asymptotic regime.

\subsection{Validity of the Superconvergent Error Estimate}
To further validate the dispersion-dissipation error estimate introduced previously, we verify that it can be used to replicate the numerical behavior of $A_T$ reported by \cite{vincent2011insights}. Namely, we use \cref{eq:final} in conjunction with \cref{eq:spectral_OA} to compute $A_T$ for some fixed $\Delta \theta$ when varying $c$. It should be noted that while the numerator of \cref{eq:final} can be evaluated analytically, its denominator must still be computed numerically since we have no closed-form expression for the eigenvalue $\lambda_1$. Hence, it can be said that this procedure yields a semi-analytical approximation of $A_T$. To perform these computations, for simplicity, $\lambda_1(\theta,c)$ was approximated by $\lambda_1(0,c)$ in \cref{eq:final}.

\begin{figure}
    \centering
    \subfloat[$p=2$ ESFR Scheme.]{
    \includegraphics[width=0.6\linewidth]{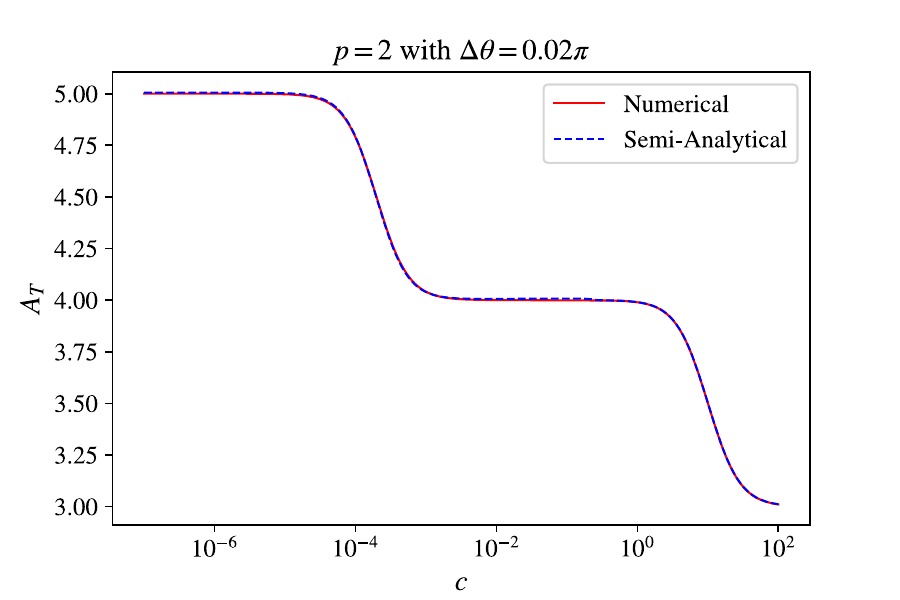}
    }
    \hfill
    \subfloat[$p=3$ ESFR Scheme.]{
    \includegraphics[width=0.6\linewidth]{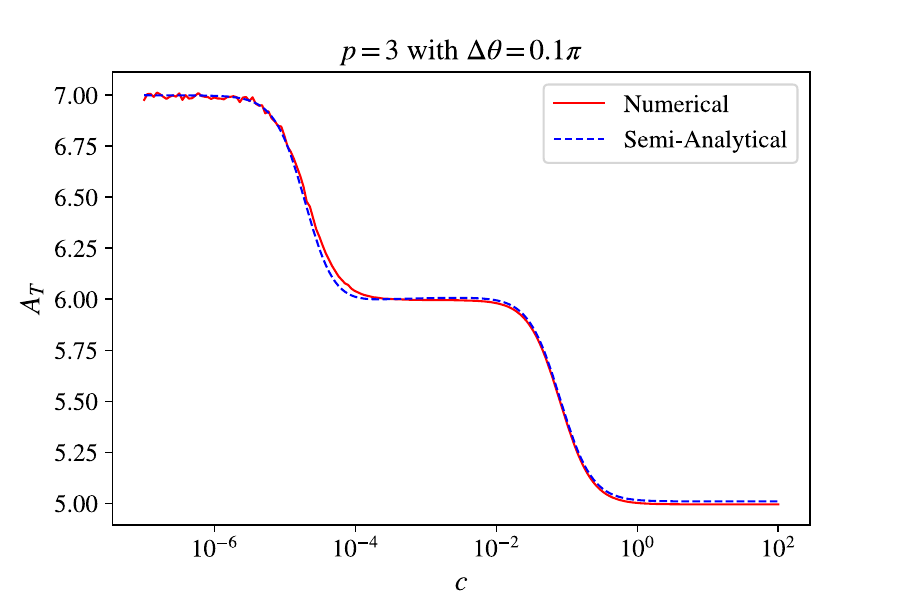}
    }
    \caption{Experimental validation of the error estimate.}
    \label{fig:ATcomparison}
\end{figure}

Results obtained for $p=2$ and $p=3$ ESFR schemes are shown in \cref{fig:ATcomparison}. For comparison, $A_T$ has also been computed purely numerically using \cref{eq:Bloch_eig} and \cref{eq:spectral_OA}, as done by \cite{vincent2011insights}. Clearly, provided that $\Delta \theta$ is sufficiently small, the error estimate given by \cref{eq:final} captures very well the numerical behavior of $A_T$. It should be noted that the values of $\Delta \theta$ used in this numerical experiment are smaller than those utilized by \cite{vincent2011insights}, as this was required to get a good agreement between the semi-analytical and purely numerical approximations of $A_T$. In both cases studied, the value of $\Delta \theta$ was set to the minimum allowable value given the applicable double machine precision limitations for the purely numerical computation of $A_T$.

\section{Extension to General Symmetric FR Schemes}
\label{sec:extension}
Although we have focused on the superconvergent properties of ESFR schemes, the methodology presented in the scope of this work can easily be extended to a very general class of FR schemes with a symmetric correction function. Indeed, as shown in \cref{app:lem_symFR}, provided that a given FR scheme has symmetric correction functions and that the projection of the derivative of the latter on each Legendre polynomial is non-vanishing, one can always find a diagonal matrix $\mathbf{K} \in \mathbb{R}^{(p+1) \times (p+1)}$ such that
\begin{equation}
    \mathbf{h}_R = (\mathbf{M} + \mathbf{K})^{-1} \mathbf{r}
    \quad ; \quad
    \mathbf{h}_L = -(\mathbf{M} + \mathbf{K})^{-1}\mathbf{l},
\end{equation}
in the Legendre basis. Using the approach from \cref{sec:proof}, the order of convergence of the dispersion-dissipation error associated with the FR schemes belonging to this class can easily be determined. Namely, we can see that the truncation error associated with the rational approximant $P(z)/Q(z)$ will lose one order for every non-zero diagonal entry of the matrix $\mathbf{K}$. Hence, we find that for this subset of symmetric FR schemes, $E_T$ is $O(\theta^{p+\bar{k}+2})$, where $(\mathbf{K})_{\bar{k}\bar{k}}$ is the non-zero diagonal entry of $\mathbf{K}$ with the smallest index.

\section{Conclusion}
In this work, we have provided a formal proof for the superconvergent behavior of the dispersion-dissipation error of ESFR schemes for the linear advection problem when using an upwind numerical flux. Using simple analytical tools, we have proven that the superconvergent properties of ESFR schemes rely on the ability of their associated von Neumann matrices to generate superconvergent rational approximants of the exponential function.

Our theoretical work has also clarified the dependency between the ESFR scalar parameter $c$ and the spectral order of accuracy $A_T$. Namely, we have shown that while the first drop in $A_T$ which is observed in numerical experiments as $c$ is increased can be attributed to a change in the structure of the superconvergent rational approximants, the second drop observed is in fact a numerical artifact stemming from the small magnitude of one of the unphysical eigenvalues in the vicinity of $\theta = 0$. Finally, the validity of the theoretical error estimate was successfully confirmed through numerical experiments.

Future work should include the study of other superconvergent properties of ESFR schemes in the context of linear and non-linear problems, as was done for the DG method in the past decades.

\section*{Acknowledgments}
The authors would like to thank the National Sciences and Engineering Research Council of Canada (NSERC).

\appendix
\section{Algebraic Multiplicity of Null Eigenvalue}
In this work, we have used rather extensively the fact that the physical eigenvalue $\omega^h$ of the von Neumann matrix has an algebraic multiplicity of one in the vicinity of $\theta = 0$. The following provides proof of this fact.
\label{app:lem_mult}
\begin{lemma}
    \label{lem:mult}
    Let $\mathbf{H} : \mathbb{R} \to \mathbb{C}^{(p+1)\times(p+1)}$ be the von Neumann matrix of an ESFR scheme. Provided that $c_{-} < c < \infty$, the algebraic multiplicity of the null eigenvalue of $\mathbf{H}(0)$ is one.
\end{lemma}
\begin{proof}
    In what follows, for concision, we will denote $\mathbf{H}(0)$ by $\mathbf{H}$. The characteristic polynomial of $\mathbf{H}$ is given by
    \begin{equation}
        \text{det}(\mathbf{H} - \lambda \mathbf{I}) = 0.
        \label{eq:char_pol}
    \end{equation}
    Thanks to the similarity in the structure of \cref{eq:modified_VonNeumann} and \cref{eq:char_pol}, we can essentially reuse the proof of \cref{lem:rational_approx} to find the characteristic polynomial of $\mathbf{H}$. Namely, removing the factor $e^{-i\theta}$ and replacing $i\theta$ by $i\lambda$ in the proof of \cref{lem:rational_approx}, one finds that the eigenvalues of $\mathbf{H}$ must satisfy
    \begin{align}
        P(i\lambda) - Q(i\lambda) =2^{-(p+1)}(i\lambda)^{p+1}
        &-\sum_{j,k=0}^p (-1)^j2^{k-p} q_j\psi^{(k)}_j(-1)(i\lambda)^{p-k} \nonumber \\
        &+\sum_{j,k=0}^p (-1)^j2^{k-p} q_j\psi^{(k)}_j(1)(i\lambda)^{p-k} = 0.
        \label{eq:char_pol1}
    \end{align}
    Reordering the summands in \cref{eq:char_pol1} yields
\begin{equation}
    2^{-(p+1)}(i\lambda)^{p+1} +
    \sum_{k=0}^pc_{k,p}(i\lambda)^{k} = 0,
\end{equation}
where
\begin{align}
    c_{k,p} :=& \sum_{j=0}^{p}(-1)^{p}2^{-k}q_{p-j}\left((-1)^{j} \psi^{(p-k)}_{p-j}(1)
        -(-1)^{j} \psi^{(p-k)}_{p-j}(-1)\right) \\
        =&
        \sum_{j=0}^{k}(-1)^{p}2^{-k}q_{p-j}\left((-1)^{j} \psi^{(p-k)}_{p-j}(1)
        -(-1)^{j} \psi^{(p-k)}_{p-j}(-1)\right).
        \label{eq:pol_coeff}
\end{align}
To prove that the null eigenvalue of $\mathbf{H}$ has an algebraic multiplicity of one, it suffices to show that $c_{0,p}=0$ and $c_{1,p} \neq 0$. Since even Legendre polynomials are symmetric and odd Legendre polynomials are anti-symmetric, it follows that
\begin{equation}
    \psi^{(m)}_{n}(1) = (-1)^{m+n}\psi^{(m)}_{n}(-1).
    \label{eq:Legendre_sym_derivative}
\end{equation}
Applying \cref{eq:Legendre_sym_derivative} to \cref{eq:pol_coeff}, one obtains
\begin{equation}
    c_{k,p} =
        \sum_{j=0}^{k}(-1)^{p}\left((-1)^{j}
        -(-1)^{k} \right)2^{-k}q_{p-j}\psi^{(p-k)}_{p-j}(1).
        \label{eq:pol_coeff_final}
\end{equation}
From \cref{eq:pol_coeff_final}, it is clear that
\begin{equation}
    c_{0,p} = (-1)^{p}\left(1
        -1 \right)2^{0}q_{p}\psi^{(p)}_{p}(1) = 0.
\end{equation}
Similarly,
\begin{equation}
    c_{1,p} = (-1)^{p}\left(1+1\right)2^{-1}q_{p}\psi^{(p-1)}_{p}(1)
    =
    \frac{(-1)^{p}(2p+1)\psi^{(p-1)}_{p}(1)}{2}f(c).
\end{equation}
Since $f(c) \neq 0$ as long as $c_{-} < c < \infty$, $c_{1,p} \neq 0$. Moreover, when $c \to \infty$, $f(c) \to 0$ and the algebraic multiplicity of the null eigenvalue is increased, as expected.
\end{proof}

\section{Symmetric FR Schemes}
\label{app:lem_symFR}
In the following, we show a simple property of symmetric FR schemes which is very useful when attempting to generalize the work presented in this paper.
\begin{lemma}
\label{lem:symFR}
    For all symmetric FR schemes of degree $p$ with correction functions such that
    \begin{equation}
        \forall_{j=0}^{p} : \int_{-1}^{1} g_R'(\xi) \psi_{j}(\xi) d\xi \neq 0
        \quad \text{and} \quad
        \int_{-1}^{1} g_L'(\xi) \psi_{j}(\xi) d\xi \neq 0,
        \label{eq:FRconstraint}
    \end{equation}
    there exists a diagonal matrix $\mathbf{K} \in \mathbb{R}^{(p+1) \times (p+1)}$ such that in the Legendre basis, the derivatives of the correction functions are given by
    \begin{equation}
    \mathbf{h}_R = (\mathbf{M} + \mathbf{K})^{-1} \mathbf{r}
    \quad ; \quad
    \mathbf{h}_L = -(\mathbf{M} + \mathbf{K})^{-1}\mathbf{l}.
    \end{equation}
\end{lemma}
\begin{proof}
    First, we recall that an FR scheme is said to be symmetric if and only if
    \begin{equation}
        g_L(\xi) = g_R(-\xi).
    \end{equation}
    As noted by \cite{vincent2015extended}, in the Legendre basis, this is equivalent to requiring
    \begin{equation}
        (\mathbf{h}_{L})_j = (-1)^{j+1} (\mathbf{h}_R)_j.
    \end{equation}
    Suppose that we are given $\mathbf{h}_R$ for a particular symmetric FR scheme satisfying \cref{eq:FRconstraint}. Then, it is clear that defining the diagonal matrix $\mathbf{K} \in \mathbb{R}^{(p+1) \times (p+1)}$ as
    \begin{equation}
        (\mathbf{K})_{jj} := 1 / (\mathbf{h_R})_j - (\mathbf{M})_{jj}
    \end{equation}
    ensures that
    \begin{equation}
        \mathbf{h}_R = (\mathbf{M} + \mathbf{K})^{-1} \mathbf{r}
        \label{eq:FRproof1}
    \end{equation}
    in the Legendre basis; \cref{eq:FRconstraint} is required so that \cref{eq:FRproof1} is well-defined. Moreover, in the Legendre basis using this choice of $\mathbf{K}$, we can see that
    \begin{equation}
        -((\mathbf{M}+\mathbf{K})^{-1}\mathbf{l})_j = -(-1)^j ((\mathbf{M} +\mathbf{K})^{-1})_j =(-1)^{j+1}(\mathbf{h}_R)_j = (\mathbf{h}_L)_j.
    \end{equation}
    Hence, we have found a diagonal matrix $\mathbf{K} \in \mathbb{R}^{(p+1)\times(p+1)}$ satisfying the required conditions.
\end{proof}

\newpage

\bibliographystyle{siamplain}
\bibliography{references}
\end{document}

%% file: ex_shared.tex
\usepackage{amssymb}
\usepackage{amsmath}
\usepackage{bm}

\usepackage{amsfonts}
\usepackage{graphicx}
\usepackage[caption=false]{subfig} 

\usepackage{epstopdf}
\usepackage{algorithmic}
\ifpdf
  \DeclareGraphicsExtensions{.eps,.pdf,.png,.jpg}
\else
  \DeclareGraphicsExtensions{.eps}
\fi

\newsiamremark{remark}{Remark}
\newsiamremark{hypothesis}{Hypothesis}
\crefname{hypothesis}{Hypothesis}{Hypotheses}
\newsiamthm{claim}{Claim}

\headers{On the Superconvergence of ESFR Schemes}{Mathias Dufresne-Piché and Siva Nadarajah}

\title{On the Superconvergence of ESFR Schemes\thanks{Submitted to the editors DATE.
\funding{This work is partially supported by the Natural Sciences and Engineering Research Council of Canada (NSERC).}}}

\author{Mathias Dufresne-Piché\thanks{Department of Mechanical Engineering, McGill University, Montréal, H3A 0G4, Canada 
  (\email{mathias.dufresne-piche@mail.mcgill.ca}, \email{siva.nadarajah@mcgill.ca}).}
\and Siva Nadarajah\footnotemark[2]}

\usepackage{amsopn}

\makeatletter
\newcommand*{\addFileDependency}[1]{
  \typeout{(#1)}
  \@addtofilelist{#1}
  \IfFileExists{#1}{}{\typeout{No file #1.}}
}
\makeatother


%% file: manuscript.bbl
\begin{thebibliography}{10}

\bibitem{adjerid2009discontinuous}
{\sc S.~Adjerid and T.~Weinhart}, {\em Discontinuous galerkin error estimation for linear symmetric hyperbolic systems}, Computer methods in applied mechanics and engineering, 198 (2009), pp.~3113--3129, \url{https://doi.org/10.1016/j.cma.2009.05.016}.

\bibitem{ahmad1998orthogonal}
{\sc F.~Ahmad}, {\em The orthogonal polynomials and the pade'approximation}, Applied Mathematics and Mechanics, 19 (1998), pp.~663--668, \url{https://doi.org/10.1007/BF02452374}.

\bibitem{ainsworth2004dispersive}
{\sc M.~Ainsworth}, {\em Dispersive and dissipative behaviour of high order discontinuous galerkin finite element methods}, Journal of Computational Physics, 198 (2004), pp.~106--130, \url{https://doi.org/10.1016/j.jcp.2004.01.004}.

\bibitem{allaneau2011connections}
{\sc Y.~Allaneau and A.~Jameson}, {\em Connections between the filtered discontinuous galerkin method and the flux reconstruction approach to high order discretizations}, Computer Methods in Applied Mechanics and Engineering, 200 (2011), pp.~3628--3636, \url{https://doi.org/10.1016/j.cma.2011.08.019}.

\bibitem{asthana2017consistency}
{\sc K.~Asthana, J.~Watkins, and A.~Jameson}, {\em On consistency and rate of convergence of flux reconstruction for time-dependent problems}, Journal of Computational Physics, 334 (2017), pp.~367--391, \url{https://doi.org/10.1016/j.jcp.2017.01.008}.

\bibitem{atkins2009super}
{\sc H.~Atkins and B.~Helenbrook}, {\em Super-convergence of discontinuous galerkin method applied to the navier-stokes equations}, in 19th AIAA Computational Fluid Dynamics, NASA, 2009, p.~3787, \url{https://doi.org/10.2514/6.2009-3787}.

\bibitem{cao2014superconvergence}
{\sc W.~Cao, Z.~Zhang, and Q.~Zou}, {\em Superconvergence of discontinuous galerkin methods for linear hyperbolic equations}, SIAM Journal on Numerical Analysis, 52 (2014), pp.~2555--2573, \url{https://doi.org/10.1137/130946873}.

\bibitem{castonguay2012new}
{\sc P.~Castonguay, P.~E. Vincent, and A.~Jameson}, {\em A new class of high-order energy stable flux reconstruction schemes for triangular elements}, Journal of Scientific Computing, 51 (2012), pp.~224--256, \url{https://doi.org/10.1007/s10915-011-9505-3}.

\bibitem{chalmers2017spatial}
{\sc N.~Chalmers and L.~Krivodonova}, {\em Spatial and modal superconvergence of the discontinuous galerkin method for linear equations}, Journal of Scientific Computing, 72 (2017), pp.~128--146, \url{https://doi.org/10.1007/s10915-016-0349-8}.

\bibitem{NIST:DLMF}
{\em {\it NIST Digital Library of Mathematical Functions}}.
\newblock \url{https://dlmf.nist.gov/}, Release 1.2.4 of 2025-03-15, \url{https://dlmf.nist.gov/}.
\newblock F.~W.~J. Olver, A.~B. {Olde Daalhuis}, D.~W. Lozier, B.~I. Schneider, R.~F. Boisvert, C.~W. Clark, B.~R. Miller, B.~V. Saunders, H.~S. Cohl, and M.~A. McClain, eds.

\bibitem{guo2013superconvergence}
{\sc W.~Guo, X.~Zhong, and J.-M. Qiu}, {\em Superconvergence of discontinuous galerkin and local discontinuous galerkin methods: eigen-structure analysis based on fourier approach}, Journal of Computational Physics, 235 (2013), pp.~458--485, \url{https://doi.org/10.1016/j.jcp.2012.10.020}.

\bibitem{hesthaven2007nodal}
{\sc J.~S. Hesthaven and T.~Warburton}, {\em Nodal discontinuous Galerkin methods: algorithms, analysis, and applications}, Springer Science \& Business Media, 2007.

\bibitem{horn2013roger}
{\sc R.~Horn}, {\em Roger and cr johnson, matrix analysis}, 2013.

\bibitem{hu2002eigensolution}
{\sc F.~Q. Hu and H.~L. Atkins}, {\em Eigensolution analysis of the discontinuous galerkin method with nonuniform grids: I. one space dimension}, Journal of Computational Physics, 182 (2002), pp.~516--545, \url{https://doi.org/10.1006/jcph.2002.7184}.

\bibitem{huynh2007flux}
{\sc H.~T. Huynh}, {\em A flux reconstruction approach to high-order schemes including discontinuous galerkin methods}, in 18th AIAA computational fluid dynamics conference, 2007, p.~4079, \url{https://doi.org/10.2514/6.2007-4079}.

\bibitem{lambert2023l2}
{\sc E.~Lambert and S.~Nadarajah}, {\em An l2-error estimate of energy stable flux reconstruction method}, 2023, \url{https://arxiv.org/abs/2309.03801}, \url{https://arxiv.org/abs/2309.03801}.

\bibitem{moura2015linear}
{\sc R.~C. Moura, S.~J. Sherwin, and J.~Peir{\'o}}, {\em Linear dispersion--diffusion analysis and its application to under-resolved turbulence simulations using discontinuous galerkin spectral/hp methods}, Journal of Computational Physics, 298 (2015), pp.~695--710, \url{https://doi.org/10.1016/j.jcp.2015.06.020}.

\bibitem{quaegebeur2019stability2}
{\sc S.~Quaegebeur and S.~Nadarajah}, {\em Stability of energy stable flux reconstruction for the diffusion problem using the interior penalty and bassi and rebay ii numerical fluxes for linear triangular elements}, Journal of Computational Physics, 380 (2019), pp.~88--118, \url{https://doi.org/10.1016/j.jcp.2018.12.017}.

\bibitem{quaegebeur2020insights}
{\sc S.~Quaegebeur and S.~Nadarajah}, {\em Insights on the coercivity of the esfr methods for elliptic problems}, Computers \& Mathematics with Applications, 80 (2020), pp.~2029--2044, \url{https://doi.org/10.1016/j.camwa.2020.09.001}.

\bibitem{quaegebeur2019stability1}
{\sc S.~Quaegebeur, S.~Nadarajah, F.~Navah, and P.~Zwanenburg}, {\em Stability of energy stable flux reconstruction for the diffusion problem using compact numerical fluxes}, SIAM Journal on Scientific Computing, 41 (2019), pp.~A643--A667, \url{https://doi.org/10.1137/18M1184916}.

\bibitem{vermeire2016properties}
{\sc B.~C. Vermeire and P.~E. Vincent}, {\em On the properties of energy stable flux reconstruction schemes for implicit large eddy simulation}, Journal of Computational Physics, 327 (2016), pp.~368--388, \url{https://doi.org/10.1016/j.jcp.2016.09.034}.

\bibitem{vincent2011insights}
{\sc P.~E. Vincent, P.~Castonguay, and A.~Jameson}, {\em Insights from von neumann analysis of high-order flux reconstruction schemes}, Journal of Computational Physics, 230 (2011), pp.~8134--8154, \url{https://doi.org/10.1016/j.jcp.2011.07.013}.

\bibitem{vincent2011new}
{\sc P.~E. Vincent, P.~Castonguay, and A.~Jameson}, {\em A new class of high-order energy stable flux reconstruction schemes}, Journal of Scientific Computing, 47 (2011), pp.~50--72, \url{https://doi.org/10.1007/s10915-010-9420-z}.

\bibitem{vincent2015extended}
{\sc P.~E. Vincent, A.~M. Farrington, F.~D. Witherden, and A.~Jameson}, {\em An extended range of stable-symmetric-conservative flux reconstruction correction functions}, Computer Methods in Applied Mechanics and Engineering, 296 (2015), pp.~248--272, \url{https://doi.org/10.1016/j.cma.2015.07.023}.

\bibitem{WANG20098161}
{\sc Z.~Wang and H.~Gao}, {\em A unifying lifting collocation penalty formulation including the discontinuous galerkin, spectral volume/difference methods for conservation laws on mixed grids}, Journal of Computational Physics, 228 (2009), pp.~8161--8186, \url{https://doi.org/https://doi.org/10.1016/j.jcp.2009.07.036}.

\bibitem{williams2013energy}
{\sc D.~M. Williams, P.~Castonguay, P.~E. Vincent, and A.~Jameson}, {\em Energy stable flux reconstruction schemes for advection--diffusion problems on triangles}, Journal of Computational Physics, 250 (2013), pp.~53--76, \url{https://doi.org/10.1016/j.jcp.2013.05.007}.

\bibitem{yang2012analysis}
{\sc Y.~Yang and C.-W. Shu}, {\em Analysis of optimal superconvergence of discontinuous galerkin method for linear hyperbolic equations}, SIAM Journal on Numerical Analysis, 50 (2012), pp.~3110--3133, \url{https://doi.org/10.1137/110857647}.

\bibitem{zwanenburg2016equivalence}
{\sc P.~Zwanenburg and S.~Nadarajah}, {\em Equivalence between the energy stable flux reconstruction and filtered discontinuous galerkin schemes}, Journal of Computational Physics, 306 (2016), pp.~343--369, \url{https://doi.org/10.1016/j.jcp.2015.11.036}.

\end{thebibliography}
